\documentclass[psamsfonts]{amsart}

\usepackage{amssymb,amsfonts}
\usepackage[all,arc]{xy}
\usepackage{enumerate}
\usepackage{mathrsfs}
\usepackage{graphicx}
\usepackage{cite}

\newtheorem{theorem}{Theorem}[section]
\newtheorem{corollary}[theorem]{Corollary}
\newtheorem{proposition}[theorem]{Proposition}

\theoremstyle{definition}
\newtheorem{definition}[theorem]{Definition}

\newtheorem{remark}[theorem]{Remark}

\DeclareMathOperator{\Nil}{Nil}

\DeclareMathOperator{\ann}{ann}


\numberwithin{equation}{section}

\begin{document}

\title{Algebraic Properties of Expectation Semirings}

\author[Peyman Nasehpour]{\bfseries Peyman Nasehpour}

\address{Peyman Nasehpour\\
	Department of Engineering Science \\
	Golpayegan University of Technology \\
	Golpayegan\\
	Iran}
\email{nasehpour@gut.ac.ir, nasehpour@gmail.com}



\begin{abstract}
In this paper, we investigate the algebraic properties of the expectation semirings which are semiring version of the concept of trivial extension in ring theory. We discuss ideals, primes, maximals and primary ideals of these semirings. We also discuss the distinguished elements such as the units, idempotents, and zero-divisors of the expectations semirings. Similar to their counterparts in ring theory, we introduce pr\'{e}simplifiable, domainlike, clean, almost clean, and weakly clean semiring and see when an expectation semiring is one of these semirings. 
\end{abstract}

\maketitle

\section{Introduction}

Semirings are ring-like algebraic structures that subtraction is either impossible or disallowed, interesting generalizations of rings and distributive lattices, and have important applications in many different branches of science and engineering. For general books on semiring theory and its applications, one may refer to the resources \cite{Golan1999(a), Golan1999(b), Golan2003, GondranMinoux2008, Glazek2002, HebischWeinert1998}.

Since different authors define semirings differently \cite{Glazek2002}, it is very important to clarify, from the beginning, what we mean by a semiring in this paper. By a semiring, we understand an algebraic structure $(S,+,\cdot,0,1)$ with the following properties:

\begin{enumerate}
	\item $(S,+,0)$ is a commutative monoid,
	\item $(S,\cdot,1)$ is a monoid with $1\neq 0$,
	\item $a(b+c) = ab+ac$ and $(b+c)a = ba+ca$ for all $a,b,c\in S$,
	\item $a\cdot 0 = 0\cdot a = 0$ for all $a\in S$. 
\end{enumerate}

A semiring $S$ is commutative if $ab = ba$ for all $a,b\in S$. All semirings except for one fleeting instance (cf. Remark \ref{ExSeRe}) will be commutative with identity. Now, let $(M,+,0)$ be a commutative additive monoid. In this paper, the monoid $M$ is said to be an $S$-semimodule if $S$ is a semiring and there is a function, called scalar product, $\lambda: S \times M \longrightarrow M$, defined by $\lambda (s,m)= sm$ such that the following conditions are satisfied:

\begin{enumerate}
	\item $s(m+n) = sm+sn$ for all $s\in S$ and $m,n \in M$;
	\item $(s+t)m = sm+tm$ and $(st)m = s(tm)$ for all $s,t\in S$ and $m\in M$;
	\item $s\cdot 0=0$ for all $s\in S$ and $0 \cdot m=0$ and $1 \cdot m=m$ for all $m\in M$.
\end{enumerate}

For more on semimodules and their subsemimodules, check Section 14 of the book \cite{Golan1999(b)}.

The expectation semiring introduced in \cite{Eisner2001} has important applications in computational linguistics and natural language processing. We borrow the general definition of such semirings (check Definition \ref{ExSeDe} in the paper) from Example 7.3 of the book \cite{Golan2003} and investigate some of the algebraic properties of these semirings. Now we clarify a brief sketch of the contents of our paper. The definition of expectation semirings is as follows:

Let $S$ be a semiring and $M$ an $S$-semimodule. On the set $S \times M$,  we define the two addition and multiplication operations as follows:

\begin{itemize}
	\item $(s_1, m_1) + (s_2,m_2) = (s_1 + s_2 , m_1 + m_2)$,
	
	\item $(s_1, m_1) \cdot (s_2, m_2) = (s_1 s_2, s_1 m_2 + s_2 m_1)$.
\end{itemize}

The semiring $S \times M$ equipped with the above operations, denoted by $S \widetilde{\oplus} M$, is called the expectation semiring of the $S$-semimodule $M$ (see Proposition \ref{ExpectationSemiring}).

Let us recall that if $S$ is a semiring and $I$ is a nonempty subset of $S$, then $I$ is called an ideal of $S$ if $a+b\in I$ and $sa\in I$, for all $a,b\in I$ and $s\in S$ \cite{Bourne1951}. Prime, maximal, and primary ideals of semirings are defined similar to their counterparts in commutative ring theory. In the same section, we investigate the prime and maximal ideals of the expectation semirings (see, for instance, Theorem \ref{TrivialExtensionThm1}).

An $S$-subsemimodule $N$ of an $S$-semimodule $M$ is subtractive if $x+y \in N$ with $x\in N$ implies that $y\in N$, for all $x,y\in M$. An $S$-semimodule $M$ is called subtractive if each subsemimodule of $M$ is subtractive. In Theorem \ref{PrimaryThm}, we investigate the primary ideals of the expectation semirings. A corollary of this theorem (see Corollay \ref{PrimaryThm2}) is as follows:

Let $S$ be a semiring, $M$ a subtractive $S$-semimodule, $I$ an ideal of $S$, and $N$ an $S$-subsemimodule of $M$. Then the ideal $I \widetilde{\oplus} N$ of the expectation semiring $S \widetilde{\oplus} M$ is primary if and only if $N$ is a primary $S$-subsemimodule of $M$, $IM \subseteq N$, and $\sqrt{I} = \sqrt{N}$, where by $\sqrt{N}$, we mean the ideal $\sqrt{[N:M]}$ of $S$.

Let us recall that a proper ideal $I$ of a semiring $S$ is defined to be weakly prime if $0 \neq ab \in I$ implies that $a\in I$ or $b\in I$ \cite{Atani2007}. In Proposition \ref{weaklyprimepro}, we show that if $S$ is a semiring, $M$ an $S$-semimodule, and $I$ a proper ideal of $S$. Then $I \widetilde{\oplus} M$ is a weakly prime ideal of the expectation semiring $S \widetilde{\oplus} M$ if and only if $I$ is weakly prime, and $ab = 0$, $a\neq 0$, and $b\neq 0$ imply that $a,b \in \ann(M)$, for all $a,b \in S$.

Noetherian semirings and semimodules are defined similar to their counterparts in module theory. In Theorem \ref{ExpectationNoetherian}, we investigate the conditions that the expectation semiring $S \widetilde{\oplus} M$ is Noetherian.

In Section \ref{sec:DistinguishedElements}, we investigate distinguished elements of the expectation semirings. For example, in Theorem \ref{TrivialExtensionThm2}, we prove the following:

Let $S$ be a semiring and $M$ an $S$-semimodule. Then the following statements hold:

\begin{enumerate}
	
	\item The set of all units $U(S \widetilde{\oplus} M)$ of $S \widetilde{\oplus} M$ is the set $U(S) \widetilde{\oplus} V(M)$, where $U(S)$ is the set of all units of the semiring $S$ and $V(M)$ is the set of all elements of the semimodule $M$ such that they have additive inverses.

	\item The set of all nilpotent elements $\Nil(S \widetilde{\oplus} M)$ of $S \widetilde{\oplus} M$ is the ideal $\Nil(S) \widetilde{\oplus} M$.
	
	\item The set of all zero-divisors $Z(S \widetilde{\oplus} M)$ of the semiring $S \widetilde{\oplus} M$ is $$\{(s,m) : s\in Z(S) \cup Z(M), m\in M \}.$$
	
\end{enumerate}

In this section, we also introduce pr\'{e}simplifiable and strongly associate semimodules: We define an $S$-semimodule $M$ pr\'{e}simplifiable if whenever $s \in S$ and $m\in M$ with $sm=m$, then either $s\in U(S)$ or $m=0$. We define a semiring $S$ to be pr\'{e}simplifiable if $S$ is pr\'{e}simplifiable as an $S$-semimodule. Then in Theorem \ref{presimplifiable}, we prove that the expectation semiring $S \widetilde{\oplus} M$ is pr\'{e}simplifiable if and only if $V(M) = M$ and $S$ and $M$ are both pr\'{e}simplifiable.

In the final phase of Section \ref{sec:DistinguishedElements}, we also introduce domainlike, clean, almost clean, and weakly clean semirings. For example, we say a semiring $S$ is clean if every element of $S$ is the sum of a unit and an idempotent (check Definition \ref{cleansemiring}). Then in Proposition \ref{cleansemiringPro}, we show that if $S$ is a semiring and $M$ an $S$-semimodule such that $V(M) = M$,then the expectation semiring $S \widetilde{\oplus} M$ is clean if and only if $S$ is clean.

We define a semiring $S$ to be almost clean if each element of the semiring may be written as the sum of a non-zero-divisor and an idempotent. Then, we prove that the expectation semiring $S \widetilde{\oplus} M$ is almost clean if and only if each element of $s\in S$ can be written in the form $s=t+e$ such that $t\notin(Z(S) \cup Z(M))$ and $e$ is an idempotent element of $S$ (see Definition \ref{almostcleansemiring} and Proposition \ref{almostcleanPro}).

Finally, we define a semiring $S$ to be weakly clean if for each $s\in S$ either $s=u+e$ or $u+e = u$, for some unit $u$ and an idempotent $e$ (check Definition \ref{weaklycleansemiring}) and we show that if $S$ is a semiring and $M$ an $S$-semimodule such that $V(M) = M$, then the expectation semiring $S \widetilde{\oplus} M$ is weakly clean if and only if $S$ is weakly clean (see Proposition \ref{weaklycleanPro}).

\section{Expectation Semirings and Their Ideals}\label{sec:expsemi}

We start this section with the following straightforward proposition:

\begin{proposition}

\label{ExpectationSemiring}

Let $S$ be a semiring and $M$ an $S$-semimodule. Then the following statements hold:

\begin{enumerate}

\item \label{trivialE} The set $S \times M$ equipped with the following addition and multiplication is a semiring: $$(s_1, m_1) + (s_2,m_2) = (s_1 + s_2 , m_1 + m_2)$$ $$(s_1, m_1) \cdot (s_2, m_2) = (s_1 s_2, s_1 m_2 + s_2 m_1).$$

\item The semiring $S$ is isomorphic to the subsemiring $S \times \{0\}$ of the semiring $S \times M$.

\item The set $\{0\} \times M$ is a nilpotent ideal of $S \times M$ and if $M \neq \{0\}$, then the index of nilpotency of $\{0\} \times M$ is 2.
    
\item \label{ExpectationMatrix} If $E$ is the set of all matrices of the form $\begin{pmatrix}
  s & m \\
  0 & s
 \end{pmatrix}$, where $s\in S$ and $m\in M$, then $E$ equipped with componentwise addition and the following multiplication $$\begin{pmatrix}
  s_1 & m_1 \\
  0 & s_1
 \end{pmatrix} \cdot \begin{pmatrix}
  s_2 & m_2 \\
  0 & s_2
 \end{pmatrix} = \begin{pmatrix}
  s_1 s_2 & s_1m_2 + s_2 m_1 \\
  0 & s_1 s_2
 \end{pmatrix},$$ is a semiring and isomorphic to the semiring defined in the statement (\ref{trivialE}) in the current proposition.

\end{enumerate}

\end{proposition}

\begin{remark}
	
	\label{trivialextension}
	
	Let $\mathbb R^{\geq 0}$ denote the set of all non-negative real numbers, $M$ be an $\mathbb R^{\geq 0}$-semimodule and define $E$ to be the set of all matrices of the form $\begin{pmatrix}
	r & m \\
	0 & r
	\end{pmatrix}$, where $r\in \mathbb R^{\geq 0}$ and $m\in M$. It is easy to verify that $E$ with the componentwise addition and the multiplication defined in the statement (\ref{ExpectationMatrix}) in Proposition \ref{ExpectationSemiring} is a called the expectation semiring by Eisner in \cite{Eisner2001}, which has important applications in computational linguistics and natural language processing.
	
	On the other hand, if $R$ is a commutative ring with a nonzero identity and $M$ a unital $R$-module, $R \times M$ equipped with the addition and multiplication defined in the statement (\ref{trivialE}) in Proposition \ref{ExpectationSemiring} is a commutative ring with the identity $(1,0)$, called a trivial extension of $R$ by $M$ \cite[p. 8]{Sernesi2006}. Apparently, the Japanese mathematician Masayoshi Nagata (1927--2008) was the first who introduced the concept of the trivial extension of a ring by a module, but under the term ``the principle of idealization'', and used that extensively in his book \cite{Nagata1962}. Note that the semiring version of the trivial extension of $R$ by $R$ has been discussed in \cite{GuptaKumar2011}, where $R$ is an arbitrary semiring.
\end{remark}

By considering Proposition \ref{ExpectationSemiring} and Remark \ref{trivialextension}, we give the following definition:

\begin{definition}

\label{ExSeDe}

Let $S$ be a semiring and $M$ be an $S$-semimodule. On the set $S \times M$,  we define the two addition and multiplication operations as follows:

\begin{enumerate}
\item $(s_1, m_1) + (s_2,m_2) = (s_1 + s_2 , m_1 + m_2)$,

\item $(s_1, m_1) \cdot (s_2, m_2) = (s_1 s_2, s_1 m_2 + s_2 m_1)$.
\end{enumerate}

The semiring $S \times M$ equipped with the above operations, denoted by $S \widetilde{\oplus} M$, is called the expectation semiring of the $S$-semimodule $M$. Note that in general if $T \subseteq S$ and $N \subseteq M$, then by $T \widetilde{\oplus} N$, we mean the set of all ordered pairs $(t,n)$ such that $t\in T$ and $n\in N$.

\end{definition}

\begin{remark}
	
	\label{ExSeRe}
	
	Let $S$ and $T$ be two semirings and $M$ be $(S,T)$-bisemimodule. Put $\begin{pmatrix}
	S & M \\
	0 & T
	\end{pmatrix}$ to be the set of all matrices of the form $\begin{pmatrix}
	s & m \\
	0 & t
	\end{pmatrix}$, where $s\in S$, $t\in T$, and $m\in M$. Define addition componentwise and multiplication as follows:
	
	$$\begin{pmatrix}
	s_1 & m_1 \\
	0 & t_1
	\end{pmatrix} \cdot \begin{pmatrix}
	s_2 & m_2 \\
	0 & t_2
	\end{pmatrix} = \begin{pmatrix}
	s_1 s_2 & s_1m_2 + t_2 m_1 \\
	0 & t_1 t_2
	\end{pmatrix}.$$
	
	It is, then, easy to verify that the set $\begin{pmatrix}
	S & M \\
	0 & T
	\end{pmatrix}$ equipped with the componentwise addition and the above multiplication is a (not necessarily commutative) semiring \cite[Example 7.3]{Golan2003}. Such semirings have applications in the automatic parallelization of linear computer codes \cite{BuckerBuschelmanHovland2002}.
	
\end{remark}

The concept of graded rings and modules is classical \cite{Bourbaki1970}. Graded semirings have been defined as almost a graded ring but with no opposite for the addition in \cite[p. 1391]{Balmer2005}. Group graded semirings \cite{SharmaJoseph2008} and $\mathbb Z$-graded semirings \cite{AllenKimNeggers2013} have been investigated recently. Similarly, we define monoid graded semirings as follows:

\begin{definition}

\label{Gradedsemiring}

Let $(M,+)$ be a commutative monoid. A semiring $S$ is said to be an $M$-graded semiring if there are submonoids $\{S_i\}_{i\in M}$ of $S$ such that $S$ is a direct sum of its submonoids $\{S_i\}_{i\in M}$ and $S_i S_j \subseteq S_{i+j}$, for all $i,j \in M$.

Also, let $I$ be an ideal of the semiring $S$ and put $I_m = I \cap S_m$, for all $m \in M$. The ideal $I$ is called $M$-graded if the following conditions hold:

\begin{enumerate}

\item $I = \bigoplus_{m\in M} I_m$ (as commutative monoids),

\item $S_k I_l \subseteq I_{k+l}$, for all $k,l \in M$.

\end{enumerate}

\end{definition}

We use the concept of $\mathbb N_0$-graded semirings in the following:

\begin{theorem}

\label{TrivialExtensionThm1}

Let $S \widetilde{\oplus} M$ be the expectation semiring of an $S$-semimodule $M$. Then the following statements hold:

\begin{enumerate}

\item The expectation semiring $S \widetilde{\oplus} M$ of an $S$-semimodule $M$ is an $\mathbb N_0$-graded semiring.

\item If $I$ is an ideal of $S$ and $N$ is an $S$-subsemimodule of $M$, then $I \widetilde{\oplus} N$ is an ideal of $S \widetilde{\oplus} M$ if and only if $IM \subseteq N$, and if it is so, then $I \widetilde{\oplus} N$ is an $\mathbb N_0$-graded ideal of $S \widetilde{\oplus} M$. On the other hand, if an ideal $J$ of $S \widetilde{\oplus} M$ is $\mathbb N_0$-graded, then there is an ideal $I$ of $S$ and an $S$-subsemimodule $N$ of $M$ such that $J = I \widetilde{\oplus} N$.

\item If $I$ is an ideal of $S$, $N$ is an $S$-subsemimodule of $M$, and $IM \subseteq N$, then $\sqrt{I \widetilde{\oplus} N} = \sqrt{I} \widetilde{\oplus} M$.

\item If $J$ is an ideal of $S \widetilde{\oplus} M$ and $I=\{s\in S: \exists ~ m\in M~(s,m)\in J\}$ and $N=\{n\in M: \exists ~ s\in S~(s,n)\in J\}$, then $I$ is an ideal of $S$, $N$ is an $S$-subsemimodule of $M$, $IM \subseteq N$, and $J \subseteq I \widetilde{\oplus} N$.

\item If $J$ is a subtractive ideal of $S \widetilde{\oplus} M$ such that includes $(0)\widetilde{\oplus} M$, then there is an ideal $I$ of $S$ such that $J = I \widetilde{\oplus} M$.

\item Each prime ideal of the semiring $S \widetilde{\oplus} M$ includes the ideal $(0) \widetilde{\oplus} M$.

\item Each subtractive prime ideal $P$ of the semiring $S \widetilde{\oplus} M$ is of the form $P=\mathfrak{p} \widetilde{\oplus} M$, where $\mathfrak{p}$ is a subtractive prime ideal of $S$.

\end{enumerate}

\begin{proof}

(1): Define $T_0 = S \widetilde{\oplus} (0)$, $T_1 = (0) \widetilde{\oplus} M$, and $T_n = ((0,0))$, for all $n\geq 2$. It is easy to see that $S \widetilde{\oplus} M = \bigoplus^{+\infty}_{i=0} T_i$ and by Proposition \ref{ExpectationSemiring}, $T_i T_j \subseteq T_{i+j}$, for all $i,j \in \mathbb N_0$.

(2): First, we prove that if $IM \subseteq N$, then $J = I \widetilde{\oplus} N$ is an ideal $S \widetilde{\oplus} M$. It is clear that $I \widetilde{\oplus} N$ is a submonoid of $S \widetilde{\oplus} M$. Now let $(s,m) \in S \widetilde{\oplus} M$ and $(a,x) \in I \widetilde{\oplus} N$. By definition, $(s,m)(a,x) = (sa, sx+am)$. Obviously $sa\in I$. Since $IM \subseteq N$, $sx+am \in N$. Now let $J = I \widetilde{\oplus} N$ be an ideal of $S \widetilde{\oplus} M$. Take $a\in I$ and $m\in M$. Definitely $(0,am) = (a,0)(0,m)\in I \widetilde{\oplus} N$. This implies that $IM \subseteq N$. In order to prove that $I \widetilde{\oplus} N$ is an $\mathbb N_0$-graded ideal of $S \widetilde{\oplus} M$, we proceed as follows: Define $J_0 = J \cap T_0 = (I \widetilde{\oplus} N) \cap (S \widetilde{\oplus} (0)) = I \widetilde{\oplus} (0),$ $J_1 = J \cap T_1 = (I \widetilde{\oplus} N) \cap ((0) \widetilde{\oplus} M) = (0) \widetilde{\oplus} N,$ and $J_n = 0$, for all $n\geq 2$. It is easy to verify that $J = \bigoplus^{+\infty}_{n=0} J_n$ as commutative monoids and $T_k J_l \subseteq J_{k+l}$, for all $k,l \in \mathbb N_0$.

Finally, let $J$ be an $\mathbb N_0$-graded ideal of $S \widetilde{\oplus} M$. Note that $T_n = 0$, for all $n\geq 2$. So $J = J_0 \oplus J_1$, where $J_n = J \cap T_n$, for $n=0,1$.

Define $I = \{s\in S: (s,0)\in J_0 \}$ and $N = \{m\in M: (0,m)\in J_1 \}$. It is, then, easy to show that $I$ is an ideal of $S$, $N$ is an $S$-subsemimodule of $M$, $J_0 = I \widetilde{\oplus} (0)$, $J_1 = (0) \widetilde{\oplus} N$ and so we have $J =  I \widetilde{\oplus} N$.

(3): Let $(s,m)\in \sqrt{I \widetilde{\oplus} N}$. Then there is a natural number $k$, such that $(s,m)^k \in I \widetilde{\oplus} N$. But $(s,m)^k = (s^k, ks^{k-1}x)$. This means that $s\in \sqrt{I}$, which implies that $(s,m) \in \sqrt{I} \widetilde{\oplus} M$. Now let $(s,m) \in \sqrt{I} \widetilde{\oplus} M$. So there is a natural number $k$ such that $s^k \in I$. Since $(s,m)^{k+1} = (s^{k+1}, (k+1)s^k m)$ and $IM \subseteq N$, $(s,m)^{k+1} \in I \widetilde{\oplus} N$, which implies that $(s,m) \in \sqrt{I \widetilde{\oplus} N}$.

(4): Let $J$ be an ideal of $S \widetilde{\oplus} M$ and put $I=\{s\in S: \exists ~ m\in M~(s,m)\in J\}$ and $N=\{n\in M: \exists ~ s\in S~(s,n)\in J\}$. It is easy to verify that $I$ is an ideal of $S$ and $N$ is an $S$-semimodule of $M$. Now let $a\in I$ and $x\in M$. It is clear that there is an $m\in M$ such that $(a,m) \in J$. Since $J$ is an ideal of $S \widetilde{\oplus} M$, we have that $(a,m)(0,x)=(0,ax)$ is an element of $J$. This implies that $ax\in N$, which means that $IM \subseteq N$. This point that $J \subseteq I \widetilde{\oplus} N$ is trivial by the definition of $I$ and $N$.

(5): Let $J$ be a subtractive ideal of $S \widetilde{\oplus} M$ and $I = \{s\in S : \exists ~ m\in M ((s,m)\in J)\}$. It is easy to check that $I$ is an ideal of $S$ and $J \subseteq I \widetilde{\oplus} M$.

Now take $(a,x)\in I \widetilde{\oplus} M$. Since $a\in I$, there is an $m\in M$ such that $(a,m)\in J$. But $J$ includes  $(0) \widetilde{\oplus} M$. So $(0,x), (0,m), (a, m+x)\in J$. Since $J$ is subtractive and $(a, m+x) = (a,x)+(0,m)$, we have $(a,x)\in J$.

(6): Let $P$ be a prime ideal of $S \widetilde{\oplus} M$. Since $((0) \widetilde{\oplus} M)^2 = 0 \subseteq P$, by primeness of $P$, we have that $(0) \widetilde{\oplus} M \subseteq P$.

(7): Let $P$ be a subtractive prime ideal of $S \widetilde{\oplus} M$. Since $P$ is prime, it contains $(0) \widetilde{\oplus} M$. Consequently, $P = \mathfrak{p} \widetilde{\oplus} M$, where $\mathfrak{p}= \{s\in S : \exists ~ m\in M ((s,m)\in P)\}$. It is, now, easy to check that $ \mathfrak{p}$ is a subtractive prime ideal of $S$.
 \end{proof}

\end{theorem}

\begin{corollary}

Let $M$ be an $S$-semimodule. Then the following statements hold:

\begin{enumerate}

\item All $\mathbb N_0$-graded ideals of the expectation semiring $S \widetilde{\oplus} M$ are subtractive if and only if the semiring $S$ and the $S$-semimodule $M$ are subtractive.

\item If $S \widetilde{\oplus} M$ is a subtractive semiring, then $S$ is a subtractive semiring and $M$ is a subtractive semimodule.

\end{enumerate}

\begin{proof}
First, we prove assertion (1). By Theorem \ref{TrivialExtensionThm1}, each $\mathbb N_0$-graded ideal of $S \widetilde{\oplus} M$ is of the form $J = I \widetilde{\oplus} N$, where $I$ is an ideal of $S$ and $N$ is an $S$-subsemimodule of $M$. Therefore, if all $\mathbb N_0$-graded ideals of $S \widetilde{\oplus} M$ are subtractive, then all ideals $I \widetilde{\oplus} (0)$ and $(0) \widetilde{\oplus} N$ are subtractive, which implies that all ideals $I$ of $S$ and all subsemimodules $N$ of $M$ are subtractive. On the other hand, if $S$ and $M$ are subtractive, then any ideal of the form $I \widetilde{\oplus} N$ is subtractive, where $I$ is an ideal of $S$ and $N$ is a subsemimodule of $M$. This obviously means that all $\mathbb N_0$-graded ideals of $S \widetilde{\oplus} M$ are subtractive. The assertion (2) is just a straightforward result of (1) and the proof is complete.
\end{proof}
\end{corollary}

Let us recall that a semiring $S$ is weak Gaussian if and only if each prime ideal of $S$ is subtractive (See Definition 18 and Theorem 19 in \cite{Nasehpour2016}).

\begin{corollary}

Let $S \widetilde{\oplus} M$ be a weak Gaussian semiring. Then the following statements hold:

\begin{enumerate}

\item Each prime ideal $P$ of the semiring $S \widetilde{\oplus} M$ is of the form $P=\mathfrak{p} \widetilde{\oplus} M$, where $\mathfrak{p}$ is a subtractive prime ideal of $S$.

 \item Each maximal ideal of the semiring $S \widetilde{\oplus} M$ is of the form $\mathfrak{m} \widetilde{\oplus} M$, where $\mathfrak{m}$ is a subtractive maximal ideal of $S$.

\end{enumerate}

\end{corollary}

\begin{theorem} 
	
\label{ExpectationNoetherian}

Let $M$ be an $S$-semimodule. Then the following statements hold:
	
	\begin{enumerate}
		
		\item If $S \widetilde{\oplus} M$ is a Noetherian semiring, then $S$ is a Noetherian semiring and $M$ is a finitely generated $S$-semimodule.
		
		\item If $S$ is Noetherian, $M$ is finitely generated, and $S \widetilde{\oplus} M$ is subtractive, then $S \widetilde{\oplus} M$ is Noetherian.
		
	\end{enumerate}
	
	\begin{proof}
		
		(1): Let $I_1 \subseteq I_2 \subseteq \cdots \subseteq I_n \subseteq$ be an ascending chain of ideals of $S$. So, $I_1 \widetilde{\oplus} (0) \subseteq I_2 \widetilde{\oplus} (0) \subseteq \cdots \subseteq I_n \widetilde{\oplus} (0) \subseteq$ is also an ascending chain of ideals of the Noetherian semiring $S \widetilde{\oplus} M$. So, there is an $i\geq 0$ such that $I_{i+n} \widetilde{\oplus} (0) = I_{i} \widetilde{\oplus} (0)$, for all $n\in \mathbb N$. So, $I_{i+n} = I_{i}$, for all $n\in \mathbb N$, which means that $S$ is Noetherian. Now, we prove that $M$ is finitely generated. Clearly, the ideal $(0) \widetilde{\oplus} M$ is finitely generated and $(0,m_1), (0,m_2), \ldots, (0,m_n)$ are its generators. Take an arbitrary element $m$ of $M$. So, $(0,m)$ is an element of $(0) \widetilde{\oplus} M$ and therefore, there are $(s_1,x_1),(s_2,x_2),\ldots,(s_n,x_n)$ in $S \widetilde{\oplus} M$ such that $$(0,m) = (s_1,x_1)(0,m_1) + (s_2,x_2)(0,m_2) + \cdots + (s_n,x_n)(0,m_n).$$ This implies that $$m = s_1 m_1 + s_2 m_2 + \cdots + s_n m_n.$$ So, $M$ is generated by $m_1,m_2,\ldots,m_n$.
		
		(2): Cohen's theorem in semiring theory states that a subtractive semiring $S$ is Noetherian if and only if every prime ideal of $S$ is finitely generated \cite[Proposition 7.17]{Golan1999(b)}. Let $P$ be a prime ideal of the subtractive semiring $S \widetilde{\oplus} M$. So, by Theorem \ref{TrivialExtensionThm1}, there is a subtractive prime ideal $\mathfrak{p}$ of $S$ such that $P = \mathfrak{p} \widetilde{\oplus} M$. On the other hand, since $S$ is Noetherian, $\mathfrak{p}$ is finitely generated. Now, since $M$ is finitely generated, $P = \mathfrak{p} \widetilde{\oplus} M$ is finitely generated and the proof is complete.
	\end{proof}
\end{theorem}

Let us recall that a nonzero nonunit element $p$ of a ring is called weakly prime if $p | ab \neq 0$ implies $p | a$ or $p | b$. Weakly prime elements have applications in factorization in rings with zero-divisors \cite{AgargunAndersonValdes1999, Galovich1978}. Based on this, a proper ideal $I$ of $R$ is defined to be weakly prime if $0 \neq ab \in I$ implies that $a\in I$ or $b\in I$. It is, then, easy to verify that the element $p$ of $R$ is weakly prime if and only if the principal ideal $(p)$ is a weakly prime ideal of $R$ \cite{AgargunAndersonValdes1999}. The semiring version of weakly prime ideals has been defined in \cite{Atani2007} as follows:

\begin{definition}
	
	\label{weaklyprimedef}
Let $S$ be a semiring. A proper ideal $I$ of $S$ is defined to be weakly prime if $0 \neq ab \in I$ implies that $a\in I$ or $b\in I$.
\end{definition}

\begin{proposition}
	
	\label{weaklyprimepro}
	
	Let $S$ be a semiring, $M$ an $S$-semimodule, and $I$ a proper ideal of $S$. Then $I \widetilde{\oplus} M$ is a weakly prime ideal of the expectation semiring $S \widetilde{\oplus} M$ if and only if $I$ is weakly prime, and $ab = 0$, $a\neq 0$, and $b\neq 0$ imply that $a,b \in \ann(M)$, for all $a,b \in S$.
	
\begin{proof}
$(\Rightarrow)$: It is easy to see that if $I \widetilde{\oplus} M$ is weakly prime, then $I$ is so. Now, suppose that $a,b\in S$ such that $ab = 0$, $a\neq 0$, and $b\neq 0$. Our claim is that $a,b \in \ann(M)$. In contrary, let $a\notin \ann(M)$. Clearly, this means that there is an $m\in M$ such that $am \neq 0$. Clearly, $(0,0) \neq (a,0)(b,m) \in  I \widetilde{\oplus} M$, while $(a,0),(b,m) \notin I \widetilde{\oplus} M$, a contradiction.

$(\Leftarrow)$: Let $a,b\in S$ and $m,n \in M$ be such that $(a,m)(b,n)\in I \widetilde{\oplus} M$. If $ab \neq 0$, then since $I$ is weakly prime, we have $(a,m)$ or $(b,m)$ is an element of $I \widetilde{\oplus} M$. Now let $ab=0$, while $a,b \notin I$. So, by assumption, $a,b \in \ann(M)$. Therefore, $(a,m)(b,n) = (0,0)$. This finishes the proof.
\end{proof}
\end{proposition}

As with module theory, one can give the following definition:

\begin{definition}

\label{PrimaryDef}

Let $M$ be an $S$-semimodule and $N$ its $S$-subsemimodule.

\begin{enumerate}

\item The residual of $M$ by $N$, denoted by $[N:M]$, is the subset $\{s\in S : sM \subseteq N\}$ of $S$.

\item The radical of $N$ in $M$, denoted by $\sqrt{N}$, is the subset $\sqrt{[N:M]}$ of $S$.

\item The subsemimodule $N$ of $M$ is called primary if $N \neq M$ and for all $s\in S$ and $m\in M$, $sm\in N$ and $m \notin N$ imply that $s^n M \subseteq N$, for some positive integer $n$.

\end{enumerate}

\end{definition}

\begin{proposition}

Let $M$ be an $S$-semimodule and $N$ its $S$-subsemimodule. Then the following statements hold:

\begin{enumerate}

\item The residual $[N:M]$ of $M$ by $N$ is an ideal of the semiring $S$.

\item If the subsemimodule $N$ of $M$ is primary, then the radical $\sqrt{N}$ of $N$ is a prime ideal of the semiring $S$.

\end{enumerate}

\begin{proof}

The first statement is straightforward. We only prove the second one: First, note that since $N \neq M$, $1 \notin \sqrt{N}$. So $\sqrt{N}$ is a proper ideal of $S$. Now let $st \in \sqrt{N}$. By definition, there is a positive integer $n$ such that $(st)^n M \subseteq N$. Suppose $t \notin \sqrt{N}$, then there is an $x\in M$ such that $t^n x \notin N$. Since $s^n t^n x \in N$, while $t^n x \notin N$, and $N$ is primary, there is a positive integer $k$ such that $s^{nk} M \subseteq N$. Thus $s\in \sqrt{N}$ and the proof is complete.
\end{proof}

\end{proposition}

\begin{theorem}

\label{PrimaryThm}

Let $S$ be a semiring, $M$ an $S$-semimodule, $I$ an ideal of $S$, and $N \neq M$ an $S$-subsemimodule of $M$. Then the following statements hold:

\begin{enumerate}

\item The ideal $I$ of $S$ is primary if and only if $I \widetilde{\oplus} M$ is a primary ideal of $S \widetilde{\oplus} M$.

\item \label{PrimarySubtractive} If the ideal $I \widetilde{\oplus} N$ of the expectation semiring $S \widetilde{\oplus} M$ is primary, then $N$ is a primary $S$-subsemimodule of $M$, $IM \subseteq N$, and $\sqrt{I} = \sqrt{N}$.

\end{enumerate}

On the other hand, if $N$ is subtractive, the converse of the statement \ref{PrimarySubtractive} holds. 

\begin{proof}

(1): For proving $\Rightarrow$, let $st\in I$ while $s \notin I$. Clearly, this implies that $(s,0)(t,0) \in I \widetilde{\oplus} N$ while $(s,0) \notin I \widetilde{\oplus} N$. Therefore, there is an $n\in \mathbb N$ such that $(s,0)^n \in I \widetilde{\oplus} N$, which means that $s^n \in I$. For the proof of $\Leftarrow$, let $(s,x)(t,y) \in I \widetilde{\oplus} M$, while $(t,y) \notin I \widetilde{\oplus} M$. Obviously, $st\in I$ and $t\notin I$, which imply that $s^n \in I$, for some $n\in \mathbb N$. So, $(s,x)^n \in I \widetilde{\oplus} M$.

(2): Now, suppose that $N \neq M$ and $sx \in N$, while $x \notin N$. Obviously, we have $(s,0)(0,x) \in I \widetilde{\oplus} N$ and $(0,x) \notin I \widetilde{\oplus} N$, which imply that $(s,0)^n \in I \widetilde{\oplus} N$, for some $n \in \mathbb N$. This means that $s^n \in I$. Note that by Theorem \ref{TrivialExtensionThm1}, we have $IM \subseteq N$, so $s^n M \subseteq N$. So, $N$ is primary.

Finally, let $x\in \sqrt{I}$. This implies that $x^n M \subseteq N$, which means that $x\in \sqrt{N}$. On the other hand, let $x\in \sqrt{N}$ and choose $m\in M-N$. Then $(x^n,0)(0,m) \in I \widetilde{\oplus} N$, for some positive integer $n$. Since $I \widetilde{\oplus} N$ is primary, $(x^n,0)^k \in I \widetilde{\oplus} N$, for some positive integer $k$, which implies that $x^{nk} \in I$.

Now, let $N$ be a subtractive and primary $S$-subsemimodule of $M$ such that $IM \subseteq N$ and $\sqrt{I} = \sqrt{N}$. We prove that $I \widetilde{\oplus} N$ is primary. Let $(s,x)(t,y) \in I \widetilde{\oplus} N$, while $(t,y) \notin I \widetilde{\oplus} N$. If $t\notin I$, then $s^n \in I$. Since $(s,x)^{n+1} = (s^{n+1}, (n+1)s^n x)$ and $IM \subseteq N$, we have $(s,x)^{n+1} \in I \widetilde{\oplus} N$. Now, if $t\in I$, then for this reason that we have supposed that $(t,y) \notin I \widetilde{\oplus} N$, we have $y\notin N$ and since $IM \subseteq N$, we have $tx\in N$. On the other hand, since $sy+tx \in N$ and $N$ is subtractive, we have $sy\in N$. By considering this point that $N$ is primary, there is a natural number $n$ such that $s^n M \subseteq N$. So once again, $(s,x)^{n+1} \in I \widetilde{\oplus} N$ and this finishes the proof.
\end{proof}

\end{theorem}

\begin{corollary}
	
	\label{PrimaryThm2}
	
	Let $S$ be a semiring, $M$ a subtractive $S$-semimodule, $I$ an ideal of $S$, and $N$ an $S$-subsemimodule of $M$. Then the ideal $I \widetilde{\oplus} N$ of the expectation semiring $S \widetilde{\oplus} M$ is primary if and only if $N$ is a primary $S$-subsemimodule of $M$, $IM \subseteq N$, and $\sqrt{I} = \sqrt{N}$.
	
\end{corollary}

\section{The Distinguished Elements of the Expectation Semirings}\label{sec:DistinguishedElements}

We recall that a nonempty subset $W$ of a semiring $S$ is called multiplicatively closed set (for short an MC-set) if $1\in W$ and for all $w_1,w_2 \in W$, we have $w_1 w_2 \in W$. An MC-set $W$ of a semiring $S$ is called saturated if $ab \in W$ if and only if $a,b\in W$ for all $a,b\in S$ \cite{Nasehpour2018P}. We also recall that if $M$ is a nonzero $S$-semimodule, an element $s\in S$ is a zero-divisor of $M$, if there is a nonzero $m\in M$ such that $sm =0$. The set of all zero-divisors of $M$ is denoted by $Z_S(M)$, or simply $Z(M)$ if there is no fear of ambiguity. 

\begin{proposition}

\label{zero-divisorunionofprimes}

Let $M$ be an $S$-semimodule. The set of zero-divisors $Z(M)$ of $M$ is a union of prime ideals of $S$.

\begin{proof}
If we take $Z(M)$ to be the set of all zero-divisors of $M$, it is easy to see that $W = S - Z(M)$ is a saturated MC-set of $S$ and by Theorem 3.5 in \cite{Nasehpour2018S}, $Z(M)$ is a union of prime ideals of $S$.
\end{proof}

\end{proposition}

Let us recall that if $(M,+)$ is a commutative monoid, then an element $x\in M$ is said to have an additive inverse if there is a $y\in M$ such that $x+y =0$. We denote the set of all elements of $M$ having additive inverses by $V(M)$. The set $V(M)$ has this property that $x+y \in V(M)$ if and only if $x\in V(M)$ and $y\in V(M)$, for all $x,y \in M$.

\begin{theorem}

\label{TrivialExtensionThm2}

Let $S$ be a semiring and $M$ an $S$-semimodule. Then the following statements hold:

\begin{enumerate}

\item The set of all units $U(S \widetilde{\oplus} M)$ of $S \widetilde{\oplus} M$ is the set $U(S) \widetilde{\oplus} V(M)$, where $U(S)$ is the set of all units of the semiring $S$ and $V(M)$ is the set of all elements of the semimodule $M$ such that they have additive inverses.

\item The element $(s,m) \in S \widetilde{\oplus} M$ is multiplicatively idempotent if and only if $s$ is multiplicatively idempotent in $S$ and $sm+sm = m$. In particular, if the only additive idempotent element of $M$ is zero, then the only multiplicatively idempotent elements of $S \widetilde{\oplus} M$ are of the form $(s,0)$, where $s$ is a multiplicatively idempotent element of $S$.

\item The set of all nilpotent elements $\Nil(S \widetilde{\oplus} M)$ of $S \widetilde{\oplus} M$ is the ideal $\Nil(S) \widetilde{\oplus} M$.

\item The set of all zero-divisors $Z(S \widetilde{\oplus} M)$ of the semiring $S \widetilde{\oplus} M$ is $$\{(s,m) : s\in Z(S) \cup Z(M), m\in M \}.$$

\end{enumerate}

\begin{proof}

(1): Let $(s,x) \in U(E)$. This means that there are $t\in S$ and $y\in M$ such that $(s,x)(t,y) = (1,0)$. This obviously implies that $st=1$ and $tx+sy =0$. It is clear that $s\in U(S)$. But $tx+sy =0$ implies that $x+s^2 y =0$, which means that $x$ has an additive inverse. On the other hand, if $(s,m)\in U(S) \times V(M)$, then $(s,m)(s^{-1}, s^{-2} (-m)) =(1,0)$.

(2): It is clear that $(s,m)^2 = (s,m)$ if and only if $s^2 = s$ and $sm+sm = m$. Now let the only additive idempotent element of $M$ be zero and $(s,m)\in S \widetilde{\oplus} M$ be idempotent. It is clear that $sm + sm = s^2m+s^2m = s(sm+sm) = sm$. This implies that $sm=0$ and so $m=0$.

(3): This is immediate from the statement (3) in Theorem \ref{TrivialExtensionThm1}.

(4): ``$\supseteq$'': First, we prove that if $s\in Z(S) \cup Z(M)$, then $(s,0) \in Z(S \widetilde{\oplus} M)$. If $s\in Z(S)$, there is a nonzero $t\in S$ such that $st=0$. Clearly $(s,0)(t,0) = (0,0)$. If $s\in Z(M)$, then there is a nonzero $m\in M$ such that $sm =0$. It is easily seen that $(s,0)(0,m) = (0,0)$ and this is the proof for what we have already claimed.

Now let $s\in Z(S) \cup Z(M)$ and $m\in M$. Since $(0,m)^2 = (0,0)$, it belongs to any prime ideal of $S \widetilde{\oplus} M$. On the other hand, since $(s,0)$ is a zero-divisor, by Proposition \ref{zero-divisorunionofprimes}, it belongs to some prime ideals of $S \widetilde{\oplus} M$ contained in $Z(S \widetilde{\oplus} M)$. From these two, we get that $(s,m) = (s,0) + (0,m)$ is an element of $Z(S \widetilde{\oplus} M)$.

``$\subseteq$'': Let $(s,m)\in Z(S \widetilde{\oplus} M)$. So for some $(t,n)\neq (0,0)$, we have $(s,m)(t,n) = (st, sn+tm) = (0,0)$. If $t\neq 0$, then $s\in Z(S)$. If $t=0$, then $n\neq 0$, while $sn = 0$, which means that $s\in Z(M)$ and the proof is complete.
\end{proof}

\end{theorem}

\begin{remark}
The first statement of Theorem \ref{TrivialExtensionThm2} is a generalization of a result for trivial extension of rings given in Lemma 1.9 in the book \cite{FossumReitenGriffith1975}.
\end{remark}

Now we apply Theorem \ref{TrivialExtensionThm2} to prove some other results similar to what we have in ring theory.

\begin{proposition}
	
	Let $k$ be a semifield and $M$ a $k$-semimodule. Then the expectation semiring $E = k \widetilde{\oplus} M$ is a local semiring.
	
	\begin{proof}
		By Theorem \ref{TrivialExtensionThm2}, $U(E) = (k-\{0\}) \times V(M)$. Now we show that $E-U(E)$ is an ideal of $E$. Let $(s,x) \in E-U(E)$. So, either $s=0$ or $x\notin V(M)$.
		
		If $s=0$, then $(t,y)(s,x) = (0,tx)$ is not a unit and therefore, is an element of $E-U(E)$.
		
		If $x\notin V(M)$ and $t \in k-\{0\}$, then $tx\notin V(M)$ and we have $tx+sy \notin V(M)$. Therefore, $(st,  tx+sy) = (t,y)(s,x)$ cannot be a unit of $E$.
		
		Finally, If $x\notin V(M)$ and $t=0$, then again $(t,y)(s,x) = (0,sy)$ is not a unit and in any case, we have that $(t,y)(s,x) \in E-U(E)$, for $(s,x) \in E-U(E)$ and $(t,y) \in E$. Therefore, by Proposition 3.5 in \cite{Nasehpour2018}, $E = k \widetilde{\oplus} M$ is local. This completes the proof.
	\end{proof}
	
\end{proposition}

Let us recall that a ring $R$ is called pr\'{e}simplifiable if whenever $r,r^{\prime} \in R$ with $rr^{\prime}=r^{\prime}$, then either $r\in U(R)$ or $r^{\prime}=0$ \cite{Bouvier1971}. An $R$-module $M$ is called pr\'{e}simplifiable if whenever $r \in R$ and $m\in M$ with $rm=m$, then either $r\in U(R)$ or $m=0$ \cite{AndersonValdesLeon1997}. Similarly, we give the following definition:

\begin{definition}
	\begin{enumerate}
		\item We define a semiring $S$ pr\'{e}simplifiable if whenever $s,t \in S$ with $st=t$, then either $s\in U(S)$ or $t=0$. 
		
		\item We define an $S$-semimodule $M$ pr\'{e}simplifiable if whenever $s \in S$ and $m\in M$ with $sm=m$, then either $s\in U(S)$ or $m=0$.
	\end{enumerate}
\end{definition}

\begin{theorem}
	
	\label{presimplifiable}
	
	Let $S$ be a semiring and $M$ an $S$-semimodule. The expectation semiring $S \widetilde{\oplus} M$ is pr\'{e}simplifiable if and only if $V(M) = M$ and $S$ and $M$ are both pr\'{e}simplifiable.
	
\begin{proof} 
$(\Leftarrow)$: Let $(s,m)$ and $(t,n)$ be elements of the expectation semiring $S \widetilde{\oplus} M$ such that $(s,m)(t,n) = (t,n)$. So, $(st, tm+sn) = (t,n)$, which means that $st = t$ and $tm+sn = n$. Since $S$ is pr\'{e}simplifiable, we have either $s\in U(S)$ or $t=0$. If $t=0$, then $sn=n$. Since $M$ is pr\'{e}simplifiable, we have either $s\in U(S)$ or $n=0$. If $n=0$, then $(t,n)=(0,0)$ and we are done. Otherwise, $s\in U(S)$. In this case, since $V(M) = M$, by Theorem \ref{TrivialExtensionThm2}, $(s,m)$ is a unit of $S \widetilde{\oplus} M$.

$(\Rightarrow)$: Let the expectation semiring $S \widetilde{\oplus} M$ be pr\'{e}simplifiable. First, we prove that $V(M)=M$. Let $m\in M$. If $m=0$, then $m$ has additive inverse and so, $m\in V(M)$. Now, let $m$ be nonzero. It is clear that $(1,m)(0,m) = (0,m)$. Since $(0,m)$ is a nonzero element of the pr\'{e}simplifiable semiring $S \widetilde{\oplus} M$, $(1,m)$ needs to be a unit of $S \widetilde{\oplus} M$. Therefore, by Theorem  \ref{TrivialExtensionThm2}, $m\in V(M)$. Finally, by considering that $st=t$ is equivalent to $(s,0)(t,0)=(t,0)$ and also, $sm=m$ is equivalent to $(s,0)(0,m)=(0,m)$, for all $s,t\in S$ and $m\in M$, it is straightforward to see that $S$ and $M$ are both pr\'{e}simplifiable and this finishes the proof. \end{proof}
		
\end{theorem}

As with module theory \cite{AndersonAxtellFormanStickles2004}, we give the following definition:

\begin{definition}
	Let $S$ be a semiring and $M$ an $S$-semimodule.
	
	\begin{enumerate}
		\item Two elements $m$ and $n$ of $M$ are associates if $(m)=(n)$, where by $(m)$, we mean the cyclic $S$-semimodule of $M$ generated by $m$. In this case, we write $m \sim n$.
		
		\item  Two elements $m$ and $n$ of $M$ are strong associates if $m=un$ for some $u\in U(S)$. In this case, we write $m\approx n$.
		
		\item We say that $M$ is strongly associate if for $m,n\in M$, $m \sim n$ implies $m\approx n$. The semiring $S$ is strongly associate if $S$ is strongly associate as an $S$-semimodule.
	\end{enumerate}
\end{definition}

\begin{proposition}
	
Each pr\'{e}simplifiable $S$-semimodule $M$ is strongly associate.
	
\begin{proof} The proof is straightforward and so, omitted.
\end{proof}
	
\end{proposition}

\begin{proposition}
	Let $S$ be a semiring and $M$ an $S$-semimodule. Then the following statements hold:
	
	\begin{enumerate}
		\item If the expectation semiring $S \widetilde{\oplus} M$ is strongly associate, then $S$ and $M$ are both strongly associate.
		
		\item Suppose that $S$ is pr\'{e}simplifiable and $V(M) = M$. Then $S \widetilde{\oplus} M$ is strongly associate if and only if $M$ is strongly associate.
	\end{enumerate}

\begin{proof} (1): Let $s$ and $t$ be elements of $S$ such that $s \sim t$. By definition, $(s) = (t)$. It is, then, easy to verify that $((s,0)) = ((t,0))$. Therefore, $(s,0) \sim (t,0)$. Since the expectation semiring $S \widetilde{\oplus} M$ is strongly associate, we have $(s,0) \approx (t,0)$. This means that there is a unit $(u,v)$ in $S \widetilde{\oplus} M$ such that $(s,0)=(u,v)(t,0)$. This implies that $s=ut$. Note that by Theorem \ref{TrivialExtensionThm2}, $u\in U(S)$. So, $s \approx t$. This means that $S$ is strongly associate. Similarly, one can prove that $M$ is also strongly associate.
	
(2): $(\Leftarrow)$: Let $M$ be strongly associate. Take two elements $(s,m)$ and $(t,n)$ in the expectation semiring $S \widetilde{\oplus} M$ such that $(s,m) \sim (t,n)$. Therefore, there two elements $(s_1,m_1)$ and $(s_2,m_2)$ in $S \widetilde{\oplus} M$ such that $(s,m)=(s_1,m_1)(t,n)$ and $(t,n) = (s_2,m_2)(s,m)$. These two imply that $(s,m) = (s_1,m_1)(s_2,m_2)(s,m)$. So, $s=s_1 s_2 s$. Now, since $S$ is pr\'{e}simplifiable, either $s_1 s_2$ is a unit of $S$ or $s=0$. If $s_1s_2$ is a unit of $S$, so $s_1$ is also a unit of $S$. On the other hand, $V(M) = M$. So, by Theorem \ref{TrivialExtensionThm2}, $(s_1,m_1)$ is a unit of $S \widetilde{\oplus} M$, which means that $(s,m) \approx (t,n)$. And if $s=0$, then $t=0$ and $m \sim n$. Now, since $M$ is strongly associate, $m \approx n$, i.e. there is a unit $u$ is $S$ such that $m=un$. So, $(0,m) = (u,0)(0,n)$, which means that $(0,m) \approx (0,n)$ and this finishes the proof.
	\end{proof}
\end{proposition}

Let us recall that  a ring $R$ is called domainlike if $Z(R) \subseteq \Nil(R)$ \cite[Definition 10]{Spellman2002}. Inspired by this, we give the following definition:

\begin{definition} Let $S$ be a semiring and $M$ an $S$-semimodule.
	
\begin{enumerate}
\item We define the semiring $S$ to be domainlike if $Z(S) \subseteq \Nil(S)$.
\item We define the $S$-semimodule $M$ domainlike if $Z(M) \subseteq \Nil(S)$.
	\end{enumerate}
 
\end{definition}

\begin{proposition}
	Let $S$ be a semiring and $M$ an $S$-semimodule. Then, the expectation semiring $S \widetilde{\oplus} M$ is domainlike if and only if $S$ and $M$ are both domainlike.
	
	\begin{proof}
	By Theorem \ref{TrivialExtensionThm2}, the proof is straightforward.	\end{proof}
\end{proposition}

Let us recall that a ring $R$ is clean if every element of $R$ is the sum of a unit and an idempotent \cite{Nicholson1977}. Similarly, we define clean semirings as follows:

\begin{definition}
	
	\label{cleansemiring}
	
We say a semiring $S$ is clean if every element of $S$ is the sum of a unit and an idempotent.
\end{definition}

\begin{proposition}
	
	\label{cleansemiringPro}
Let $S$ be a semiring and $M$ an $S$-semimodule such that $V(M) = M$. The expectation semiring $S \widetilde{\oplus} M$ is clean if and only if $S$ is clean.
	
\begin{proof} 
The implication $\Rightarrow$ is obvious. Now, let $(s,m) \in S \widetilde{\oplus} M$. Since by assumption $S$ is clean, $s\in S$ is the sum of a unit $u\in U(S)$ and an idempotent element $e\in S$. Clearly, $(s,m) = (u,m)+(e,0)$. By Theorem \ref{TrivialExtensionThm2}, $(u,m)$ is a unit and $(e,0)$ is an idempotent and this finishes the proof.
\end{proof}
\end{proposition}

A ring $R$ is called almost clean if each element of $R$ may be written as the sum of a regular element (an element which is not a zero-divisor) and an idempotent \cite[Definition 11]{McGovern2003}. Similarly, we define almost clean semirings as follows:

\begin{definition}
	
	\label{almostcleansemiring}
	
	We define a semiring $S$ to be almost clean if each element of the semiring may be written as the sum of a non-zero-divisor and an idempotent.
\end{definition}

\begin{proposition}
	
	\label{almostcleanPro}
Let $S$ be a semiring and $M$ an $S$-semimodule. The expectation semiring $S \widetilde{\oplus} M$ is almost clean if and only if each element of $s\in S$ can be written in the form $s=t+e$ such that $t\notin(Z(S) \cup Z(M))$ and $e$ is an idempotent element of $S$.

\begin{proof}
	
$(\Rightarrow)$: Let $S \widetilde{\oplus} M$ be almost clean. Suppose $s\in S$. Clearly, for $(s,0)\in S \widetilde{\oplus} M$, there are two elements $(t,m)$ and  $(e,n)$ in $S \widetilde{\oplus} M$ such that $(t,m)$ is non-zero-divisor and $(e,n)$ is idempotent and $(s,0) = (t,m) + (e,n)$. Using Theorem \ref{TrivialExtensionThm2}, we get that $t\notin(Z(S) \cup Z(M))$ and $e$ is an idempotent element of $S$.

$(\Leftarrow)$: Now, let $S$ be almost clean. Suppose $s\in S$ and $m\in M$. By assumption, $s=t+e$ such that $t\notin(Z(S) \cup Z(M))$ and $e$ is idempotent. On the other hand, $(s,m) = (t,m) + (e,0)$. Clearly, $(t,m)$ is non-zero-divisor and $(e,0)$ is idempotent and the proof is complete.
	\end{proof}
\end{proposition}

Let us recall that a ring $R$ is weakly clean if for each $r\in R$ either $r=u+e$ or $r=u-e$, for some unit $u$ and an idempotent $e$. Inspired by this, we give the following definition:

\begin{definition}
	
	\label{weaklycleansemiring}
	We say a semiring $S$ is weakly clean if for each $s\in S$ either $s=u+e$ or $u+e = u$, for some unit $u$ and an idempotent $e$.
\end{definition}

\begin{proposition}
	
	\label{weaklycleanPro}
Let $S$ be a semiring and $M$ an $S$-semimodule such that $V(M) = M$. The expectation semiring $S \widetilde{\oplus} M$ is weakly clean if and only if $S$ is weakly clean.

\begin{proof} The proof is similar to the proof of Proposition \ref{cleansemiringPro} and omitted.
\end{proof}
\end{proposition}

Let us recall that an element $a$ of a semiring $S$ is additively regular if there exists
an element $b$ of $S$ satisfying $a + a + b = a$ \cite[\S 13]{Golan1999(b)}. Similarly, one can define an element $x$ of an $S$-semimodule $M$ to be additively regular if there exists
an element $y$ of $M$ satisfying $x + x + y = x$. The semiring $S$ ($S$-semimodule $M$) is additively regular if each element of $S$ ($M$) is additively regular. It is straightforward to see the following:

\begin{proposition}
Let $M$ be an $S$-semimodule. Then the following statements hold:

\begin{enumerate}
	
	\item An element $(a,m)$ of the expectation semiring $S \widetilde{\oplus} M$ is additively regular if and only if $a$ is an additively regular element of $S$ and $m$ an additively regular element of $M$.
	
	\item The expectation semiring $S \widetilde{\oplus} M$ is additively regular if and only if $S$ and $M$ are both additively regular. 
	
	\end{enumerate}
	
\end{proposition}

\section*{Acknowledgments} The author is supported by the Department of Engineering Science at the Golpayegan University of Technology and his special thanks go to the Department for providing all the necessary facilities available to him for successfully conducting this research.

\bibliographystyle{plain}

\end{document}